\documentclass[12pt]{article}

\usepackage{amsfonts,amsmath,amssymb}

\numberwithin{equation}{section}

\title{"Nice" Rational Functions }

\author{Allan J. MacLeod\\Statistics, O.R. and Mathematics Group,\\School of Science,\\
University of the West of Scotland,\\High St.,  Paisley,\\Scotland.  PA1 2BE\\
(e-mail: allan.macleod@uws.ac.uk) }

\date{}

\begin{document}

\maketitle

\begin{abstract}
{\noindent We consider simple rational functions
$R_{mn}(x)=P_m(x)/Q_n(x)$, with $P_m$ and $Q_n$ polynomials of
degree $m$ and $n$ respectively. We look for "nice" functions,
which we define to be ones where as many as possible of the roots, poles, critical
points and (possibly) points of inflexion are integer or, at worst, rational. }

\vspace{0.5cm}

\end{abstract}

\section{Introduction}
The concept of "nice" polynomials has been around for several years.
These are polynomials with integer roots, critical points and, possibly,
points of inflexion, which are useful in providing elementary Calculus problems.
Several references can be found in the bibliography, with a very nice summary
given in Buchholz and MacDougall \cite{newdio}.

One stage more complex than a polynomial is a rational function,
which gives the student practice in the quotient rule and, often, the product rule.
The possible function graphing possibilities are also much greater, due to
 the consideration of poles and asymptotes, both horizontal and vertical.

In general, we consider functions of the form
\begin{equation}
R_{mn}(x)=\frac{P_m(x)}{Q_n(x)}
\end{equation}
where $m$ and $n$ denote the degrees of the relevant polynomials. We
assume true rational functions so $m,n \ge 1$, with
$P_m$ and $Q_n$ having no common factors, and both monomials, so that any rational roots or poles are
integers. We also assume that, when $m,n>1$, there are no repeated roots.

We would like any zeros and poles to be small integers, together with the
critical points, where $R'(x)=0$, and, if possible, the points of
inflexion, where $R''(x)=0$. We, further, assume that $m+n \ge 3$ to ensure
at least one critical point. We also note that $Q_n/P_m$ gives the same equation
for critical points as $P_m/Q_n$, but this does not happen for the second derivatives.

We make extensive use of Cardano's method for solving a cubic, which we now summarize. All
cubics can, by linear transformation, be put in the form
\begin{displaymath}
z^3+Qz+R=0
\end{displaymath}
and, defining $z=u+v$, $v=-Q/(3u)$ we have
\begin{displaymath}
u=\left( \frac{-R}{2} \pm \sqrt{\left( \frac{R}{2}\right)^2+\left(\frac{Q}{3}\right)^3} \right)^{1/3}
\end{displaymath}
and we denote the term under the square-root sign as $D$. We have, if $D>0$ there is only one real root,
whilst if $D<0$ there are $3$.

\section{$m+n=3$}

\subsection{$P_{21}$}
We first consider $R_{21}=P_2/Q_1$. We assume the quadratic either has $2$ integer roots or has
$2$ complex roots but integer coefficients. We can, by transformation, assume $Q_1=x$ so firstly let
\begin{equation}
R_{21} = \frac{(x-a)(x-b)}{x}
\end{equation}
with $a,b \in \mathbb Z$. We could consider complex roots by allowing $a,b$ to be complex conjugates,
but we prefer to keep things simple and integral as far as possible.

We have $R'=(x^2-ab)/x^2$ and $R''=2ab/x^3$, so we have no points of inflexion, and integer roots
iff $ab=\Box$, thus $a,b$ are both positive or both negative.

So for example $a=1,b=4$, and putting $x=z-n, n \in \mathbb Z$, gives the parametric function
\begin{displaymath}
\frac{z^2-(2n+5)z+n^2+5n+4}{z-n}
\end{displaymath}
which has roots at $z=n+1,n+4$, a single pole at $z=n$ and $2$ critical points at $z=n \pm 2$.

\vspace{1cm}

The other possibility for $R_{21}$ is
\begin{equation}
R_{21} = \frac{x^2+cx+d}{x}
\end{equation}
where $c,d \in \mathbb Z$ and $c^2<4d$.

We have $R'=(x^2-d)/x^2$ and $R''=2d/x^3$, and since $d>0$ there are always $2$ real critical points.
We need $d=\Box$ for integer values, and again there are no points of inflexion.

So, for example, setting $c=2,d=4$ gives
\begin{displaymath}
\frac{z^2+(2-2n)z+n^2-2n+4}{z-n}
\end{displaymath}
which has complex roots at $z=n-1 \pm\sqrt{3}i$, a single pole at $z=n$ and $2$ critical points at $z=n\pm2$.

\subsection{$R_{12}$}
For $R_{12}$ we assume $P_1=x$, so let
\begin{equation}
R_{12} = \frac{x}{(x-a)(x-b)}
\end{equation}

As we stated before, the condition for $R'=0$ is the same as for equation (2.1), namely $ab=\Box$, but
\begin{equation}
R''=\frac{ 2(x^3  - 3abx + ab(a + b))}{(x - a)^3(x - b)^3}
\end{equation}
so we must solve $x^3  - 3abx + ab(a + b)=0$ for points of inflexion. We have
\begin{displaymath}
D=\frac{a^2b^2(a-b)^2}{4}
\end{displaymath}
which implies that the cubic has only one real root since $a,b \ne 0$ and $a\ne b$. We get integer $u$ by having
$a^2b$ or $ab^2$ a cube. So (possibly swapping $a,b$) we have $2$ conditions for $a,b$, namely
\begin{displaymath}
ab=f^2 \hspace{2cm} a^2b=g^3
\end{displaymath}
and, it is not hard to see that we need $a=\alpha s^6, b=\alpha t^6$.

The simplest example is $a=1,b=64$ gives
\begin{equation}
R_{12} = \frac{z-n}{z^2-(2n+65)z+n^2+65n+64}
\end{equation}
which has a zero at $z=n$, $2$ poles at $z=n+1,n+64$, $2$ critical points at $z=n \pm 8$ and
a single point of inflexion at $z=n-20$.

\vspace{1cm}

Finally, we consider
\begin{equation}
R_{12}=\frac{x}{x^2+cx+d}
\end{equation}
with $c^2<4d$, which has
\begin{displaymath}
R_{12}'=\frac{d-x^2}{(x^2+cx+d)^2} \hspace{2cm} R_{12}''=\frac{2(x^3-3dx-cd)}{(x^2+cx+d)^3}
\end{displaymath}

Thus, integer critical points requires $d=e^2$ which we now assume. The cubic for $R_{12}''=0$, has
$D=d^2(c^2-4d)/4<0$, so there are always $3$ real points of inflexion. Following Cardano,
\begin{equation}
u^3=\left( \frac{cd}{2} \pm \frac{d}{2} \, \sqrt{4d-c^2} \, \imath \right)
\end{equation}
and, expressing in polar form,
\begin{displaymath}
u^3=e^3(\cos A \pm \imath \sin A)
\end{displaymath}
where $\tan A=\sqrt{4d-c^2}/c$.

Thus
\begin{displaymath}
u+v=2e\cos(A/3)
\end{displaymath}
is a solution of $R_{12}''=0$. Thus we need $\cos(A/3) \in \mathbb Q$. However, to find $\cos(A/3)$
from $\cos A$ requires solving a cubic - a nice circular argument!

Since $4d-c^2=4e^2-c^2$, it might be thought that Pythagorean triples might be useful. Numerical
tests, however, show that $d=e^2=k^6, k \in \mathbb{Z}$ give the only solutions, and we have
been unable to determine a formula for the corresponding $c$ values. The simplest are
$(c,d)=(9,64),(11,64)$, the first of which gives the function
\begin{equation}
R_{12}=\frac{z-n}{z^2+(9-2n)z+n^2-9n+64}
\end{equation}
which has a single root at $z=n$, no poles, two critical points at $z=n \pm 8$, and one integer
point of inflexion at $z=n-12$.

If $c=286, d=7^6$, there are $3$ integer points of inflexion.

\section{$m+n=4$}

\subsection{$R_{22}$}
Since both numerator and denominator are quadratic, we have $4$ possible mixtures of real roots and
complex roots, ignoring the possibility of repeated roots. Thus, in general,
\begin{equation}
R_{22}=\frac{(x-p)(x-q)}{(x-r)(x-s)}
\end{equation}
where $(p,q)$ and $(r,s)$ are either distinct integer pairs or complex-conjugate pairs giving a quadratic
with integer coefficients.

Thus
\begin{displaymath}
R_{22}'=\frac{(p+q-r-s)x^2+2(rs-pq)x+p(q(r+s)-rs)-qrs}{(x-r)^2(x-s)^2}
\end{displaymath}
and $R_{22}''=2T_3/((s-x)^3(x-r)^3)$ where
\begin{displaymath}
T_3=(p+q-r-s)x^3+3x^2(rs-pq)x^2+3(p(q(r+s)-rs)-qrs)x+
\end{displaymath}
\begin{displaymath}
p(rs(r+s)-q(r^2+rs+s^2))+rs(q(r+s)-rs)
\end{displaymath}
so we usually assume $p+q \ne r+s$.

For critical points, the discriminant of $R_{22}'=0$ is
\begin{equation}
(p-r)(p-s)(q-r)(q-s).
\end{equation}
It is easy algebra to show that if we have complex conjugate numbers then this is positive, so there are $2$ real critical points.
The only chance for $0$ critical points is certain integer values, such as $(p,q,r,s)=(2,4,1,3)$.

For points of inflexion, we can show that the cubic $T_3$ gives
\begin{equation}
D=\frac{(p-r)^2(p-s)^2(q-r)^2(q-s)^2(r-s)^2}{4(p+q-r-s)^4}
\end{equation}
and, again, it is standard algebra to show $D>0$ when $(r,s)\in \mathbb Z$ whilst $D<0$ if $(r,s)$ complex. Thus there
are $3$ real points of inflexion when the denominator of $R_{22}$ has no real roots, but only $1$ real point of inflexion if it
has real roots.

We now consider specific combinations of roots. Firstly, assume
\begin{equation}
R_{22}=\frac{(x-a)(x-b)}{x(x-c)}
\end{equation}
with $a,b,c \in \mathbb{Z}$ and $a \ne b \ne c \ne a$.

Thus
\begin{displaymath}
R_{22}'=\frac{(a+b-c)x^2-2abx+abc}{x^2(x-c)^2}
\end{displaymath}
\begin{displaymath}
R_{22}''=\frac{2((a+b-c)x^3-3abx^2+3abcx-abc^2)}{x^3(x-c)^3}
\end{displaymath}

The discriminant of the numerator of $R_{22}'$ must be an integer square for rational roots, so
\begin{displaymath}
a b(c-a)(c-b)=\Box
\end{displaymath}
which has many solutions.

For $R_{22}''=0$
\begin{equation}
D=\frac{a^2b^2c^2(c-a)^2(c-b)^2}{4(a+b-c)^4}
\end{equation}
which is, obviously, positive, so the cubic has only one real root. To have rational roots we need
\begin{displaymath}
a^2b^2(c-a)(c-b)=f^3
\end{displaymath}
with $f \in \mathbb{Z}$.

There are several "small" solutions to the two conditions on $a,b,c$. We look for small values of $|a+b-c|$,
and the smallest is $a=1, b=5, c=21$. This gives the general rational function
\begin{equation}
\frac{z^2-(2n+6)z+n^2+6n+5}{z^2-(2n+21)z+n^2+21n}
\end{equation}
with zeros at $z=n+1,n+5$, poles at $z=n,n+21$, critical points at $z=n-3,n+7/3$, and one real point of
inflexion at $z=n-7$.

\vspace{1cm}

Next, assume
\begin{equation}
R_{22}=\frac{x(x-a)}{x^2+cx+d}
\end{equation}
where $c^2-4d<0$, so no poles exist.

Thus,
\begin{displaymath}
R_{22}'=\frac{(a+c)x^2+2dx-ad}{(x^2+cx+d)^2} , \hspace{0.2cm}
R_{22}''=\frac{-2((a+c)x^3 + 3dx^2  - 3adx - d(ac +d))}{ (x^2  + c·x + d)^3}
\end{displaymath}

For rational critical points, we require
\begin{displaymath}
d(a^2+ac+d)=\Box=e^2
\end{displaymath}
whilst, for $R_{22}''=0$, we have
\begin{displaymath}
D=\frac{ d^2 (a^2  + ac + d)^2 (c^2  - 4d)}{ 4(a + c)^4}
\end{displaymath}
so, clearly, $D<0$ and hence there are $3$ real points of inflexion.

As with the final part of the previous section, we can show that one of the roots of the cubic is
given by
\begin{displaymath}
\frac{2e\cos(A/3)}{a+c}
\end{displaymath}
where $\cos A=(ac+c^2-2d)/2e$.

We use a simple search procedure and find several $(a,c,d)$ sets which give the required conditions. The simplest
is $a=3, c=2, d=5$, giving the general form
\begin{equation}
R_{22}=\frac{z^2-(2n+3)z+n^2+3n}{z^2+(2-2n)z+n^2-2n+5}
\end{equation}
which has zeros at $z=n,n+3$, no poles, two critical points at $z=n-3,n+1$ and one integer point of inflexion at
$z=n-1$.

During our search, we found some $(a,c,d)$ sets which gave $3$ rational points of inflexion. They are given in Table $3.1$.
It should be noted that the factorizations from these values might well test weaker students.

\begin{center}TABLE $3.1$ \\$(a,c,d)$ sets with $3$ pts. of inflexion\\[0.5em]
\begin{tabular}{rrr}
$a \,$ & $c \,$ & $d \,$ \\
49&-4&196\\
169&-1&169\\
196&32&931\\
343&17&343\\
539&-19&637\\
560&10&133
\end{tabular}
\end{center}

\vspace{1cm}

We now consider
\begin{equation}
R_{22}=\frac{x^2+cx+d}{x(x-a)}
\end{equation}
where $c^2<4d$. The first derivative gives the same condition for rational critical points as the previous example, namely
\begin{displaymath}
d(a^2+ac+d)=\Box=e^2
\end{displaymath}
but the second derivative is different, with points of inflexion coming from
\begin{displaymath}
(a+c)x^3+3dx^2-3adx+a^2d=0
\end{displaymath}

We have
\begin{displaymath}
D=\frac{a^2d^2(a^2+a c+d)^2}{4(a+c)^4}
\end{displaymath}
so, since $D>0$, there is only one real point of inflexion. Using Cardano's method gives
\begin{displaymath}
u^3=\frac{-d^2(a^2+ac+d)}{a+c)^3}
\end{displaymath}
so, for the single real point of inflexion to be rational, we require $d^2(a^2+ac+d)=f^3$.

The smallest solution is $a=21,c=3,d=8$ leading to
\begin{equation}
R_{22}=\frac{z^2+(3-2n)z+n^2-3n+8}{z^2-(2n+21)z+n^2+21n}
\end{equation}
which has no zeros, $2$ poles at $z=n,n+21$, two critical points at $z=n-3,n+7/3$ and one point of
inflexion at $z=n-7$.

\vspace{1cm}

The final $m=2,n=2$ form is
\begin{equation}
R_{22}=\frac{x^2+ax+b}{x^2+cx+d}
\end{equation}
with $a^2<4b \, , \, c^2<4d$, so there are no zeros or poles.

\begin{displaymath}
R_{22}'=-\frac{(a - c)x^2 + 2(b - d)x - ad + bc}{(x^2  + c·x + d)^2}
\end{displaymath}
so the condition for rational critical points is
\begin{displaymath}
(a - c)(ad - bc) + (b - d)^2=\Box
\end{displaymath}

Since,
\begin{displaymath}
R_{22}''=\frac{ 2((a - c)x^3 + 3(b - d)x^2 + 3(b c - a d)x - acd + b(c^2  - d) + d^2 )}{(x^2  + cx + d)^3}
\end{displaymath}
we find
\begin{displaymath}
D=\frac{ (c^2  - 4d)((a - c)(ad - bc) + (b - d)^2 )^2}{ 4(a - c)^4}
\end{displaymath}
so that $D<0$ implying that there are $3$ real points of inflexion.

A computer search finds many $(a,b,c,d)$ sets giving $2$ rational critical points and $1$ rational
point of inflexion. The simplest solution is $a=-1, b=2, c=2, d=5$ leading to the general form
\begin{equation}
R_{22}=\frac{z^2-(2n+1)z+n^2+n+2}{z^2+(2-2n)z+n^2-2n+5}
\end{equation}
which has $2$ critical points when $z=n-3,n+1$ and a rational point of inflexion at $z=n-1$.

The search for $(a,b,c,d)$ giving $3$ rational points of inflexion gave no solutions for small sums of the
parameter sizes. This initially led to the idea that such solutions might not exist. Refining the search
algorithm, however, finally found a plethora of solutions, which clearly had simple sequential forms.

We found the following three functions
\begin{displaymath}
R_{22}=\frac{x^2+2mx+m^2+26}{x^2+2(m-15)x+m^2-30m+333}
\end{displaymath}
\begin{displaymath}
R_{22}=\frac{x^2+(2n+1)x+n^2+n+168}{x^2+(2n-23)x+n^2-23n+301}
\end{displaymath}
\begin{displaymath}
R_{22}=\frac{x^2+2px+p^2+76}{x^2+2(p-9)x+p^2-18p+381}
\end{displaymath}
where $m,n,p \in \mathbb{Z}$, and we suspect there are an infinite number of such formulae.

\subsection{$R_{31}$}
Since $R_{31}$ has a cubic numerator there is always at least one real root, which we assume is an integer. The
quadratic remainder is assumed to have $2$ integer roots or $2$ complex roots.

Let,
\begin{equation}
R_{31}=\frac{(x-a)(x-b)(x-c)}{x}
\end{equation}
with $a,b,c \in \mathbb{Z}$.

Thus,
\begin{displaymath}
R_{31}'=\frac{2x^3-(a+b+c)x^2+abc}{x^2} \hspace{2cm} R_{31}''=\frac{2(x^3-abc)}{x^3}
\end{displaymath}
so that integer points of inflexion only occur if $abc=f^3$, and there will only be one real point of inflexion.

For $R_{31}'=0$ the cubic has
\begin{displaymath}
D=\frac{abc ( 27abc - (a+b+c)^3)}{2^4 \, 3^3}
\end{displaymath}
so we have $3$ real critical points if $(a+b+c)^3/(abc)>27$, which is certainly the case when all $3$ values are positive.

Searches have not been able to find $(a,b,c)$ giving $3$ rational critical points and $1$ rational point of inflexion. We have
$3$ critical points for $(10,15,24)$ giving the general parametric form
\begin{equation}
R_{31}=\frac{z^3-(3n+49)z^2+(3n^2+98n+750)z-n^3-49n^2-750n-3600}{z-n}
\end{equation}
with roots at $z=n+10,n+15,n+24$, one pole at $z=n$, and $3$ critical points at $z=n+12,n+20,n-15/2$.

The triple $(3,18,-32)$ gives the form
\begin{equation}
R_{31}=\frac{z^3+(11-3n)z^2+(3n^2-22n-618)z-n^3+11n^2+618n+1728}{z-n}
\end{equation}
with roots at $z=n-32,n+3,n+18$, one pole at $z=n$, one real critical value at $z=n+8$ and one real pt. of inflexion
at $z=n-12$.

The other possibility for $R_{31}$ is
\begin{equation}
R_{31}=\frac{(x-a)(x^2+cx+d)}{x}
\end{equation}
where $a,c,d,\in \mathbb Z$ and $c^2<4d$.

The equations for critical points and points of inflexion are, respectively,
\begin{displaymath}
2x^3+(c-a)x^2+ad=0 \hspace{2cm} x^3-ad=0
\end{displaymath}

Searching quickly finds many examples with $3$ rational critical points, for example the parametric form
\begin{equation}
\frac{z^3-(3n+13)z^2+(3n^2+26n+15)z-(n^3+13n^2+15n+36)}{z-n}
\end{equation}
which has $1$ integer zero at $z=n+12$, one pole at $z=n$ and $3$ rational critical points at
$z=n+2,n+6,n-3/2$.

\section{$m+n=5$}

We concentrate mainly on $R_{32}$, since $R_{32}'=0$ involves a quartic and $R_{32}''=0$ a cubic. $R_{23}$ has respectively
a quartic and a sextic, $R_{14}$ has a quartic and a septic, and $R_{41}$ has two quartics.

Firstly, assume
\begin{equation}
R_{32}=\frac{x^3+ax^2+bx+c}{(x-d)(x-e)}
\end{equation}
with $a,b,c,d,e \in \mathbb Z$,
then the equation for points of inflexion is a cubic which can be shown to have
\begin{displaymath}
D=\frac{(d-e)^2(ad^2+bd+c+d^3)^2(ae^2+be+c+e^3)^2}{4(a(d+e)+b+d^2+de+e^2)^4}
\end{displaymath}
so, if $d,e$ are real, $D>0$ and there is only ever one real point of inflexion with $2$
real poles.

We, now, consider the situation where both the numerator
and denominator have $3$ and $2$ integer roots respectively, so the general form is
\begin{equation}
R_{32}=\frac{(x-a)(x-b)(x-c)}{x(x-d)}
\end{equation}
where $\{a,b,c,d\} \in \mathbb Z$ with $a<b<c$, and none of the roots are the same as the poles.

Finding the critical points means solving
\begin{equation}
x^4-2dx^3-(a(b+c-d)+b(c-d)-c d)x^2+2a b c x - a b c d=0
\end{equation}
whilst the single real point of inflexion satisfies
\begin{equation}
(a(b+c-d)+(b-d)(c-d))x^3-3a b c x^2 + 3 a b c d x - a b c d^2 = 0
\end{equation}

The quartic for $R_{32}'=0$ means that we could have $4$,$2$ or $0$ real critical points. An
investigation of the underlying geometry show that we are guaranteed $4$ real critical points when
\begin{itemize}
\item $0<d<a<b<c$
\item $0<a<b<c<d$
\item $a<b<c<d<0$
\end{itemize}

The search region is greatly reduced using these extra conditions. Searching for $4$ rational critical
points has so far been unsuccessful. We have found several examples with $2$ integer critical points, the simplest
giving rise to the parametric form
\begin{equation}
R_{32}=\frac{z^3-(3n+13)z^2+(3n^2+26n+46)z-n^3-13n^2-46n-48}{z^2-(2n+18)z+n^2+18n}
\end{equation}
which has roots at $z=n+2,n+3,n+8$, poles at $z=n,n+18$ and $2$ integer critical points at
$z=n-2,n+6$.

For $4$ rational critical points, we need to consider the general form
\begin{equation}
w=\frac{z^3+rz^2+sz+t}{(z-d)(z-e)}
\end{equation}
where $d,e,r,s,t, \in \mathbb Z, d \ne e$. We use a variation of a method proposed by Dave Rusin, then of Northern Illinois
University, in a submission to the sci.math newsgroup \cite{rusin}.

Firstly, by using the transformations
\begin{displaymath}
z=\frac{(d-e)x+d+e}{2} \hspace{2cm} w=\frac{(d-e)y+3d+3e+2r}{2}
\end{displaymath}
we can transform the rational function to the simpler form
\begin{equation}
y=\frac{x^3+bx+c}{x^2-1}
\end{equation}
with $b,c \in \mathbb Q$. To prevent reducibility, we assume $b+c+1\ne0$ and $c-b-1\ne0$.

Thus, critical points satisfy
\begin{displaymath}
x^4-(b+3)x^2-2cx-b=0
\end{displaymath}
and suppose there are $3$ distinct rational solutions $x_1,x_2,x_3$ - the fourth root $=-b/(x_1x_2x_3)$ is automatically rational.
Using $x_1$ and $x_2$ we can eliminate $b,c$, and substitute into the third equation. Express $x_1=Y/W, x_2=Z/W, x_3=X/W$, then
\begin{equation}
3W^4-(X^2+Y^2+Z^2+XY+XZ+YZ)W^2+XYZ(X+Y+Z)=0
\end{equation}

We prodeed by first solving for $W^2$ and then trying to find $W^2=\Box$. For rational $W^2$ we must have
\begin{displaymath}
(X^2+Y^2+Z^2+XY+XZ+YZ)^2-12XYZ(X+Y+Z)=\Box
\end{displaymath}
which can be written as
\begin{displaymath}
T^2=X^4+2(Y+Z)X^3+(3Y^2-8YZ+3Z^2)X^2+
\end{displaymath}
\begin{displaymath}
2(Y+Z)(Y^2-5YZ+Z^2)X+(Y^2+YZ+Z^2)^2
\end{displaymath}

This quartic in $X$ clearly has solutions, so is birationally equivalent to an elliptic curve, see Silverman and Tate \cite{siltate}.
We find the curve to be
\begin{equation}
V^2=U(U+3(Y-Z)^3)(U+(3Y+Z)(Y+3Z))
\end{equation}
with the reverse transformation
\begin{equation}
X=\frac{V-(Y+Z)U-3(Y+Z)(Y-Z)^2}{2(U+3(Y-Z)^2)}
\end{equation}

For integer values of $Y,Z$ this curve has clearly $3$ finite points of order $2$ when
$U=0,-3(Y-Z)^2,-(3Y+Z)(Y+3Z)$. Numerical tests indicate these are the only torsion points apart
from some special cases. These torsion points lead to either undefined $X$ or values of $X$ which
give $b+c+1=0$ or $c-b-1=0$, which we previously forbad. We thus need points of infinite order.

Such points are easily found, for example
\begin{enumerate}
\item $((Y-Z)^2,\pm4(Y+Z)(Y-Z)^2)$,
\item $(-3(Y+Z)^2,\pm12YZ(Y+Z))$,
\item $(3(Y+3Z)(3Y+Z),\pm12(Y+Z)(Y+3Z)(3Y+Z))$.
\end{enumerate}

As an example, taking the last point (with the + sign) gives $X=(Y^2+4YZ+Z^2)/(Y+Z)$ and the quartic in $W$ factors
\begin{displaymath}
(W^2-(Y^2+3YZ+Z^2)\;(3(Y+Z)^2W^2-2YZ(Y^2+4YZ+Z^2))=0
\end{displaymath}

The first factor gives the quadratic form $W^2=Y^2+3YZ+Z^2$, which can be parameterized using standard methods. After clearing
denominators, we find
\begin{eqnarray*}
X=p^4-8p^3q+14p^2q^2-4pq^3-2q^4\\
Y=(3q^2-2pq)(p^2-2pq+2q^2)\\
Z=(p^2-q^2)(p^2-2pq+2q^2)\\
W=(p^2-3pq+q^2)(p^2-2pq+2q^2)
\end{eqnarray*}

A very interesting example comes from $p=4,q=3$, which gives, after some simplifying transformations,
\begin{equation}
R_{32}=\frac{x^3+165x^2+10566x+476280}{x^2+110x}
\end{equation}
which has one real zero at $x=-108$, two poles at $x=0,-110$, four critical points at $x=-126, -90, -70, 66$, and
one real point of inflexion at approximately $x=-81.461$.

Further examples of this type from this parameterisation have not been found. We, thus, resort to, for fixed $Y,Z \in \mathbb Z$,
searching the elliptic curve for as many independent points as can be found. We, actually, search the corresponding
minimal curve for integer coordinate points and transform. Suppose we find points $P_1,P_2,\ldots,P_m$. We then compute
\begin{displaymath}
Q=(U,V)=n_1P_1+\ldots+n_mP_m+T_k
\end{displaymath}
where $-3\le n_i \le 3$ and $T_k$ is one of the torsion points, including the point at infinity. From $Q$, we generate $X$ and then
see if a rational $W$ can be found. If so, we compute $b,c$ and look at the rational function. We search for situations where
either $R_{32}=0$ at a rational point or $R_{32}''=0$ at a rational point.

These are very rare, but we found the previous example quickly. The first example for a rational point of inflexion comes from
$Y=41, Z=13$, which gives, again after simplifying transformations
\begin{equation}
R_{32}=\frac{x^3+77x^2+292018x-3891096}{x^2+154x}
\end{equation}
which has one real root at approximately $x=13.27$, two poles at $x=0,-154$, four critical points at $x=-714,-34,66,374$
and one real point of inflexion at $x=2618/23$.

Other examples similar to the last two can be found. The "Holy Grail" would be an $R_{32}$ with $3$ rational zeros,
$2$ rational poles, $4$ rational critical points and $1$ rational point of inflexion. We conjecture this cannot occur, but
have no idea how to prove it!

For $R_{23}$, we can consider the form
\begin{displaymath}
y=\frac{x^2-1}{x^3+bx+c}
\end{displaymath}
which gives the same equation for critical points as $(4.7)$. The equation for points of inflexion is, however,
\begin{displaymath}
x^6-3(b+2)x^4-7cx^3-3bx^2+3cx+c^2-b^2=0
\end{displaymath}

We have been unable to use the previous elliptic curve approach, with this form, to find an $R_{23}$ with one rational point of inflexion
and $4$ rational critical points.

\end{document}